\documentclass[%
reprint,
superscriptaddress,
 amsmath,amssymb,
 aps,
 prl,
floatfix,
]{revtex4-1}

\usepackage[english]{babel}
\usepackage{graphicx}
\usepackage{dcolumn}
\usepackage{bm}
\usepackage{physics,comment}
\usepackage{xcolor}
\usepackage{float}
\usepackage{csquotes}
\usepackage[hidelinks]{hyperref}
\usepackage{cleveref}

\begin{document}

\title{Saddle avoidance of noise-induced transitions in multiscale systems}

\author{Reyk B\"{o}rner}
\thanks{R.B. and R.D. contributed equally to this work.}
\affiliation{Department of Mathematics and Statistics, University of Reading, Reading, UK}

\author{Ryan Deeley}
\thanks{R.B. and R.D. contributed equally to this work.}
\affiliation{Theoretische Physik/Komplexe Systeme, ICBM, Carl von Ossietzky Universität Oldenburg, Oldenburg, Germany}

\author{Raphael R\"{o}mer}
\affiliation{Department of Mathematics and Statistics, University of Exeter, Exeter, UK}

\author{Tobias Grafke}
\affiliation{Mathematics Institute, University of Warwick, Coventry, UK}

\author{Valerio Lucarini}
\affiliation{School of Computing and Mathematical Sciences, University of Leicester, Leicester, UK}

\author{Ulrike Feudel}
\affiliation{Theoretische Physik/Komplexe Systeme, ICBM, Carl von Ossietzky Universität Oldenburg, Oldenburg, Germany}

\date{\today}

\begin{abstract}
In multistable dynamical systems driven by weak Gaussian noise, transitions between competing states are often assumed to pass via a saddle on the separating basin boundary. By contrast, we show that timescale separation can cause saddle avoidance in non-gradient systems. Using toy models from neuroscience and ecology, we study cases where sample transitions deviate strongly from the instanton predicted by Freidlin-Wentzell theory, even for weak finite noise. We attribute this to a flat quasipotential and present an approach based on the Onsager-Machlup action to aptly predict transition paths.
\end{abstract}


\maketitle

Multistable systems, when randomly perturbed, may undergo transitions between their coexisting attracting states \cite{feudel_complex_2008, pisarchik_control_2014}. Examples range from brain activity \cite{brinkman_metastable_2022} and gene regulation \cite{bressloff_stochastic_2017} to lasers \cite{masoller_noise-induced_2002},  planetary atmospheres \cite{bouchet_rare_2019} and the earth system \cite{boers_theoretical_2022, rousseau_punctuated_2023}. Transitions often represent high-impact low-probability events, e.g. financial crashes \cite{bouchaud_langevin_1998}, ecosystem collapse \cite{bashkirtseva_sensitivity_2011}, and climate tipping points \cite{ ashwin_tipping_2012,lucarini_transitions_2019, ditlevsen_tipping_2010}.  Understanding critical transitions is crucial to assess a system's stability and resilience \cite{holling_resilience_1973}. 

Noise-driven systems are commonly formulated as
stochastic differential equations (SDEs) of It\^{o} type,
\begin{align} \label{eq:system-sde}
    \dd \bm{x} = \bm{b}(\bm{x}) \, \dd t + \sigma \bm \Sigma(\bm{x})\, \dd \bm W_t \,, \quad \bm x(0) = \bm x_0 \,, \ t \geq 0 \ ,
\end{align}
where $\bm x(t) \in \mathbb{R}^D$ evolves under the combined effect of the deterministic drift $\bm b(\bm x(t)): \mathbb{R}^D \to \mathbb{R}^D$ and a stochastic forcing by a $D$-dimensional Wiener process $\bm W_t$ scaled with noise amplitude $\sigma>0$. We assume that the  system $\dot{ \bm x} = \bm b(\bm x)$ possesses multiple attractors (hence being multistable) and consider non-degenerate noise, ensuring the covariance matrix $\bm Q(\bm{x}) = \bm \Sigma(\bm{x}) \bm\Sigma^\top(\bm{x}) \in \mathbb R^{D\times D}$ is invertible \cite{lucarini_global_2020}.
This modeling framework is fundamental to climate physics (Hasselmann's programme \cite{hasselmann_stochastic_1976, arnold_hasselmanns_2001,lucarini_theoretical_2023}), chemical physics \cite{dykman_large_1994}, theoretical biology \cite{bressloff_stochastic_2017, holling_resilience_1973}, and further applications of complex physics \cite{masoller_noise-induced_2002, de_souza_noise-induced_2007, schafer_escape_2017, dykman_theory_1979}.

The theory and intuition about systems described by Eq. \eqref{eq:system-sde} is often guided by the \textit{gradient} case, where $\bm b(\bm x) = -\nabla V(\bm x)$, and isotropic noise is assumed. The potential $V : \mathbb{R}^D \to \mathbb{R}$ then describes an energy landscape that visualizes the system's global stability, with attractors located at local minima of the landscape and \textit{saddles} marking \enquote{mountain passes} for noise-induced transitions -- where trajectories are most likely to cross between different valleys. However, most systems of interest are out of equilibrium, meaning no potential $V$ exists whose negative gradient is $\bm b$. 

For this general \textit{non-gradient} case, Freidlin-Wentzell (FW) theory \cite{freidlin_random_1998} introduces the \textit{quasipotential} as an energy-like scalar field \cite{freidlin_random_1998, cameron_finding_2012, zhou_quasi-potential_2012, zhou_construction_2016, graham_nonequilibrium_1986}. Following a large deviation principle \cite{touchette_large_2009}, the quasipotential is computed via a variational approach of minimizing the FW \textit{action} functional $S_T$ (see below). Intuitively, this functional measures the \enquote{energetic cost} of a trajectory in the limit of vanishing noise. As $\sigma \to 0$, transitions between competing attractors concentrate around a minimum action path, or \textit{instanton}, and the probability of deviating from the instanton decays exponentially \cite{freidlin_random_1998, grafke_numerical_2019}. 

In most cases, the instanton from one attractor to another passes through a saddle of $\bm b$ where the quasipotential has a minimum along the boundary separating the different basins of attraction \cite{margazoglou_dynamical_2021, lucarini_global_2020, lucarini_edge_2017}. Transition rates can be computed in terms of the quasipotential value at the saddle \cite{bouchet_generalisation_2016}. These results of FW theory have established the wide-spread view of saddles acting as gateways of noise-induced transitions.

However, the $\sigma \to 0$ limit is never attained in reality; various counter-examples attest that if the noise is weak yet non-infinitesimal, noise-induced transitions do not necessarily pass near a saddle  \cite{maier_escape_1993, maier_transition-rate_1992, luchinsky_observation_1999, berezhkovskii_rate_1989, northrup_saddlepoint_1983, agmon_dynamics_1987, schafer_escape_2017}. Here we show that saddle avoidance can occur when non-gradient dynamics features fast and slow degrees of freedom. We demonstrate that sample transition paths may deviate significantly from the FW instanton even for noise so weak that transitions become extremely rare. We resolve the apparent disagreement between FW theory and observations using the quasipotential and introducing a variational formulation based on the Onsager-Machlup (OM) action \cite{stratonovich_probability_1971, durr_onsager-machlup_1978, horsthemke_onsager-machlup_1975, pinski_transition_2010, gladrow_experimental_2021, li_gamma-limit_2021}.

\paragraph{Predicting Noise-Induced Transitions.}
Rare events tend to be predictable in the sense that they will most probably occur in the least unlikely way. In our context, noise-induced transitions between attractors of $\dot {\bm x} = \bm b(\bm x)$ become increasingly rare as $\sigma\to 0$. FW theory quantifies this via a large deviation principle: consider the set $\mathcal{C}_{IF}^T$ of paths leading from a state $\bm x_I$ to state $\bm x_F$ in time $T$. For $\sigma\to 0$, the probability that a solution $\bm \phi_t$ to Eq. \eqref{eq:system-sde} initialized at $\bm x_I$ remains inside a $\delta$-tube of a path $\bm \varphi_t\in \mathcal{C}_{IF}^T$ follows
\begin{align} \label{eq:large-deviation-principle}
    \mathbb{P}_\sigma\bigg(\sup_{0\leq t\leq T} ||\bm \phi_t-\bm \varphi_t||<\delta\bigg) \overset{\sigma\downarrow 0}{\asymp} \exp\left(- \frac{S_T[\bm \varphi_t]+\varepsilon_\delta}{\sigma^2} \right) \ ,
\end{align}
where $\asymp$ denotes asymptotic logarithmic equivalence and $\varepsilon_\delta$ satisfies $\lim_{\delta\downarrow 0}\varepsilon_\delta = 0$ \cite{freidlin_random_1998, dembo_large_2010}. The rate function $S_T[\bm \varphi_t]$ is the FW \textit{action functional},
\begin{align} \label{eq:fw-action}
  S_T[\bm\varphi_t] = \frac{1}{2} \int_0^T  || \dot {\bm\varphi}_t - \bm b(\bm\varphi_t) ||^2_{\bm Q(\bm\varphi_t)} \dd t \,,
\end{align} 
which measures the ``work'' done against $\bm b$, as weighted via the  $\bm Q$-metric $\lvert\lvert \bm v \rvert\rvert_{\bm Q} := \sqrt{\langle \bm v, \bm Q^{-1} \bm v \rangle}$.
This implies that in the limit $\sigma\to 0$, the most probable path connecting $\bm x_I$ to $\bm x_F$ is by far the minimum action path, or \textit{instanton}, $\bm\varphi^{* IF} := \arg \min_{T>0,\,\bm\varphi_t \in \mathcal{C}_{IF}^T} S_T[\bm\varphi_t]$ {\cite{grafke_numerical_2019}}. Starting from a reference attractor $\bm x_R = \bm x_I$, the ``difficulty'' of reaching any state $\bm x_F$ is given by the FW \textit{quasipotential},
\begin{align} \label{eq:fw-quasipotential}
    V_{R}(\bm x_F) = \inf_{T>0} \ \inf_{\bm\varphi_t \in \mathcal{C}_{RF}^T} S_T[\bm\varphi_t] \ .
\end{align} 
By definition, $V_R(\bm x_R)=0$ is the only local and thus global minimum.

Now suppose Eq. \eqref{eq:system-sde} has two stable equilibria $\bm x_R$, $\bm x_L$ with basins of attraction $\mathcal{B}_R$, $\mathcal{B}_L \subset \mathbb{R}^D$, separated by the basin boundary $\partial \mathcal{B}$. We define a transition path $\bm\phi_t^{RL}$ as a trajectory that, after exiting a small neighborhood $R$ around $\bm x_R$, enters a small neighborhood $L$ around $\bm x_L$ without re-entering $R$. Transitions can reach $\bm x_L$ at zero action after crossing the boundary $\partial \mathcal B$ (by following a flow line of $\bm b$); therefore, as $\sigma\to 0$ the most probable location to cross $\partial \mathcal B$ is at the global minimum of $V_R$ when restricting to $\partial \mathcal B$. This minimum is typically a saddle of $\bm b$, and determines the Kramers-like scaling law of the \textit{mean first-exit time} $\expval{\tau_\sigma^{RL}}$, i.e. the expected waiting time until a trajectory initialized at $\bm x_R$ leaves $\mathcal{B}_R$, through $ \expval{\tau_\sigma^{RL}} \asymp \exp\left(\,\sigma^{-2}\min_{\boldsymbol{x}\in\partial\mathcal{B}} V_{R}(\boldsymbol{x}) \right)$ \cite{freidlin_random_1998, arrhenius_uber_1889, kramers_brownian_1940, berglund_kramers_2013, day_large_1990, dembo_large_2010}. This emphasizes the relevance of saddles for noise-induced transitions as $\sigma\to 0$.

Apart from the mean first-exit time, a second timescale characterizing the transition is the \textit{mean transition time} $\expval{t^{RL}_\sigma}$, i.e. the average time $\bm\phi_t^{RL}$ takes to travel to $L$ after last leaving $R$. Large values of $r_\sigma^{RL} := \expval{\tau_\sigma^{RL}}/\expval{t^{RL}_\sigma}$ indicate that individual transitions are \textit{rare} and therefore occur as a memory-less Poisson process. 
In what follows we choose $\sigma$ sufficiently small such that $r_\sigma^{RL} \gg 1$, in which case one might expect FW theory to apply. We investigate the transition behavior in two paradigmatic two-dimensional multiscale models, one driven by additive and another by multiplicative noise.

\begin{figure*}
    \centering
    \includegraphics[width=\textwidth]{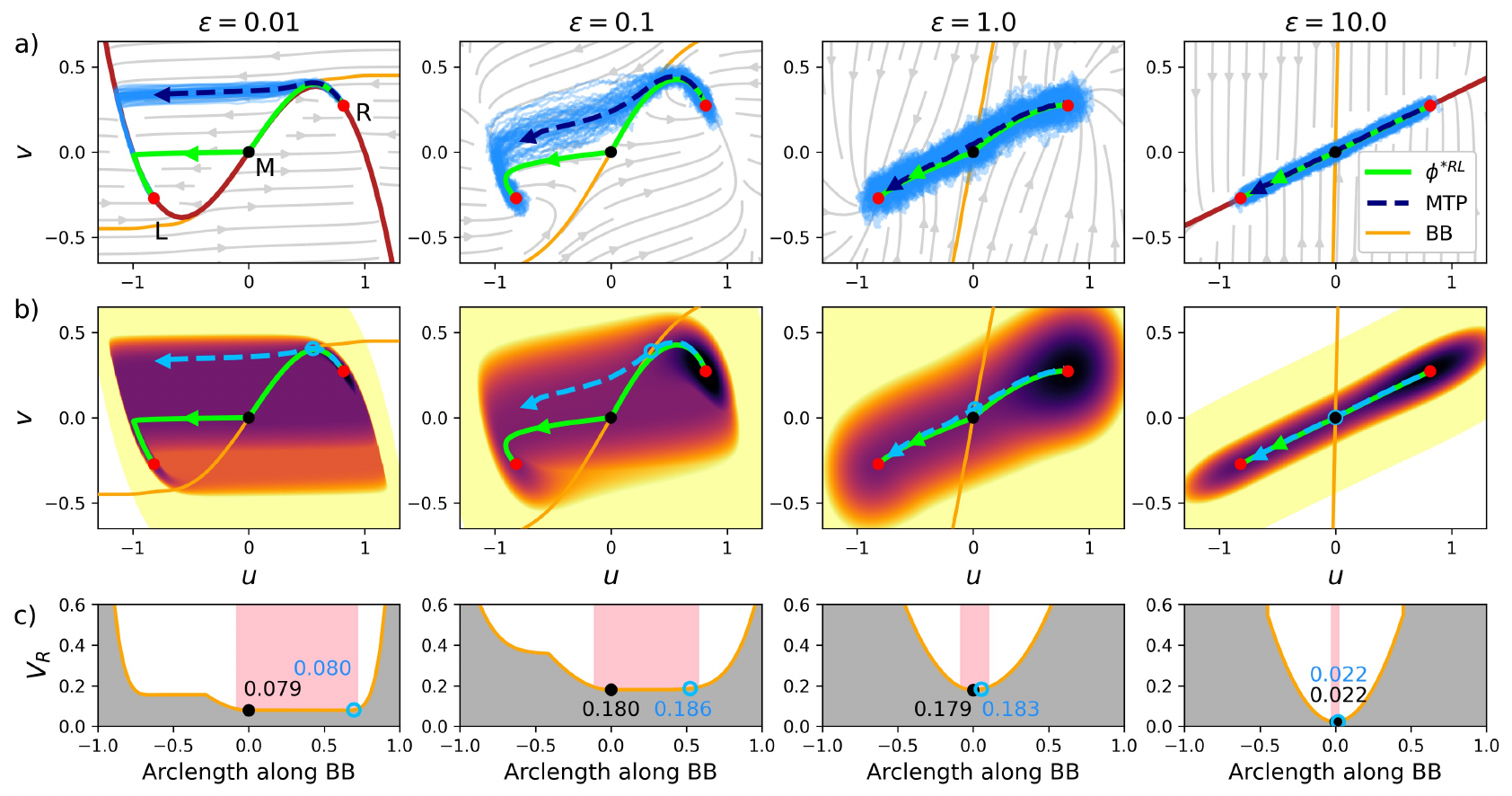}
    \caption{\label{fig:wide}FHN model for $\varepsilon = (0.01, 0.1, 1, 10 )$ at $\sigma(\varepsilon) \approx (0.08, 0.12, 0.12, 0.04)$. a) Phase space with equilibria $\bm x_{L,R}$ (red), saddle $\bm x_M$ (black), drift $\bm b$ (gray flow lines), instanton $\bm\varphi^{*RL}$ (green), 50 sample transition paths (light blue), MTP (blue dashed), and basin boundary (BB, orange); critical manifold $C_0$ (red) for $\varepsilon\to 0$ (left) and $1/\varepsilon\to 0$ (right panel). b) Like a) but with $V_R(\bm x)$ on logarithmic colormap (black: $V_R=0$; brighter $=$ larger $V_R$), showing $\bm x_c$ (blue point) where MTP crosses BB. c) $V_R(\bm x)$ along BB, indicating values at $\bm x_M$ (black) and $\bm x_c$ (blue), and the region $\mathcal{Z}$ (red shading).}
    \label{fig:fhn-main}
\end{figure*}

\paragraph{FitzHugh-Nagumo Model.} Let us first consider the FitzHugh-Nagumo (FHN) model, which was originally conceived to describe a spiking neuron \cite{fitzhugh_impulses_1961, nagumo_active_1962} and has been widely studied as a multiscale conceptual model in theoretical neuroscience. The model can be written as:
\begin{align} \label{eq:b-fhn}
    \dot{\bm{x}}=\bm b_{\text{FHN}}(u,v) = \begin{pmatrix}
        \varepsilon^{-1} (-u^3 + u - v) \\    
        -\beta v + u
    \end{pmatrix}\,.
\end{align}
Here $\bm x = (u,v) \in \mathbb{R}^2$ is the two-dimensional state vector, $\varepsilon > 0$ denotes the timescale parameter, and we set $\beta = 3$. We consider Eq. \eqref{eq:system-sde} with $\bm b = \bm b_\text{FHN}$ and identity covariance matrix $\bm Q = \bm{I}_2$ (additive noise). The noise-free system is bistable, possessing stable equilibrium points at $\bm x_{R,L} = \pm (\sqrt{2/3}, \sqrt{2/3^3})$ and a saddle point at $\bm x_M = (0,0)$. 

For $\varepsilon \ll 1$, Eq. \eqref{eq:b-fhn} describes a \textit{fast-slow} system where $u$ is \textit{fast}, while $v$ is \textit{slow} \cite{kuehn_multiple_2015}. In the $u$-direction, deterministic trajectories rapidly approach the \textit{critical manifold} $C_0 := \left\{(u,v) \in \mathbb R^2: v = - u^3 + u \right\}$ to which the slow dynamics is confined as $\varepsilon \to 0$ \cite{kuehn_multiple_2015}.
The cubic form of $C_0$ yields two saddle-node bifurcations with respect to $v$. Two stable branches of $C_0$ are separated by an unstable branch between the fold points $v_f^{\pm} = \pm \sqrt{4/3^3}$.

For such stochastic fast-slow systems with $0 < \varepsilon \ll 1$, one anticipates sample trajectories to closely track the stable part of $C_0$ until they approach a fold point, where they may abruptly transition to the opposite stable branch \cite{kuehn_mathematical_2011, zhang_predictability_2015}. By contrast, FW theory predicts that sample transition paths pass arbitrarily near the saddle for sufficiently small $\sigma$ \cite{freidlin_random_1998}. We thus face two competing limits, $\sigma \to 0$ and $\varepsilon\to 0$, when studying weakly noise-driven multiscale systems \cite{kuehn_general_2022}.

Using Monte-Carlo simulations, we sample 100 transitions $\bm\phi_t^{RL}$ for each $\varepsilon \in \{ 0.01,\, 0.1,\, 1,\, 10\}$, fixing $\sigma(\varepsilon)$ to maintain $r_{\sigma}^{RL} \approx 10^5$ (details in Supplementary Material (SI) \footnote{Supplemental Material at \url{https://figshare.com/s/281b6e7f7b556f5bbea6}.}). The ensembles of transition paths concentrate within a tube around a \textit{mean sample transition path} (MTP), which we compute by spatial averaging over the ensemble (Fig. \ref{fig:fhn-main}a). Additionally, using the geometric minimum action method (gMAM) \cite{heymann_geometric_2008}, we compute the corresponding instantons $\bm\varphi^{* RL}$ which minimize the FW action   functional $S_T$ (Eq. \eqref{eq:fw-action}).

The computed instantons always pass through the saddle. They approach it at an angle relative to $\partial \mathcal B$ that narrows as $\varepsilon$ decreases, becoming tangential to $\partial \mathcal B$ when $\varepsilon < 1/\beta$ \cite{maier_limiting_1997} (Fig. \ref{fig:fhn-main}a). 
For $\varepsilon \in \{1, 10 \}$, the MTP closely matches the instanton and all sample paths cross the basin boundary within $\sim  \sigma$ distance on either side of the saddle. 
Contrarily, for $\varepsilon \in \{0.01, 0.1\}$ sample transitions avoid the saddle: the MTP diverges from the instanton after getting close to the basin boundary \cite{maier_transition-rate_1992}, which happens before reaching the saddle. Once the noise kicks the trajectory into the competing basin, it is repelled from the basin boundary stronger than it is attracted towards the saddle. 

This multiple timescale effect manifests itself in the ratio of the negative and positive eigenvalues $\lambda_\mp$ of the Jacobian of $\bm b$ at the saddle, $\mu := |\lambda_-|/\lambda_+$. If $\mu < 1$, Ref. \cite{maier_limiting_1997} showed that sample transition paths avoid the saddle on an extended lengthscale $\mathcal{O}(\sigma^\mu)$ \cite{bobrovsky_results_1992, day_cycling_1994}: their exit locations (where they cross the basin boundary) follow a one-sided Weibull distribution \cite{coles2001introduction}, whose mode approaches the saddle only logarithmically as $\sigma \to 0$ \cite{maier_limiting_1997, luchinsky_observation_1999}. This contrasts with the case $\mu > 1$ where the distribution of exit locations is centered around the saddle on a lengthscale $\mathcal{O}(\sigma)$ as $\sigma\to 0$ \cite{maier_limiting_1997}. 
In the FHN model, $\mu < 1$ if and only if $0 < \varepsilon < 1/\beta$, and we have $\mu \approx 2\varepsilon$ for $\varepsilon \ll 1$, which highlights the link between saddle avoidance and timescale separation (see Fig. S1 in SI).
The case $\varepsilon=0.01$ exemplifies the fast-slow behavior anticipated for $\varepsilon\to 0$:
the MTP escapes from the basin boundary near the bifurcation point $v_f^+$ where $C_0$ becomes unstable (Fig. \ref{fig:fhn-main}a). In the inverse limit $1/\varepsilon \to 0$ (where $\mu \gg 1$), the corresponding critical manifold given by $\beta v=u$ is globally stable; the limiting behavior of both $\sigma$ and $1/\varepsilon$ aligns by forcing the instanton and sample transition paths onto $C_0$, as seen for $\varepsilon=10$.

Can we understand the observed transition behavior by means of the \textit{quasipotential} $V_R$ (Eq. \eqref{eq:fw-quasipotential}), which measures the “difficulty” of reaching a point $\bm x$ from $ \bm x_R$? As $\sigma\to 0$, the probability of passing through the global minimum  $\bm x^*$ of $V_R$ along the basin boundary $\partial \mathcal{B}$ approaches 1 \cite{freidlin_random_1998}. For finite $\sigma$, however, the large deviation principle underlying FW theory (Eq. \eqref{eq:large-deviation-principle}) suggests that transitions may cross with  similar probability in regions $\mathcal{Z}(\sigma) = \{\bm z \in \partial \mathcal B : V_R(\bm z) \leq V_R(\bm x^*) + \sigma^2\}$.

We compute $V_R$ for the four $\varepsilon$-cases considered using the OLIM4VAD algorithm \cite{dahiya_ordered_2018} (see SI). For $\varepsilon=10$, a steep and narrow trench of low quasipotential connects $\bm x_R$ with $\bm x_L$ (Fig. \ref{fig:fhn-main}b). As $\varepsilon$ decreases, however, $V_R$ increasingly flattens along the $v$-direction, leading to an extended \textit{quasipotential plateau} around the saddle $\bm x_M$. 
Along $\partial\mathcal{B}$, $V_R$ indeed assumes a global minimum at $\bm x_M$ for all $\varepsilon$, but its curvature around $\bm x_M$ becomes small for $\varepsilon \ll 1$ (Fig. \ref{fig:fhn-main}c). The low curvature implies the lack of a clear \enquote{mountain pass} and a widened \enquote{danger zone} $\mathcal{Z}(\sigma)$ around $\bm x_M$, whose arclength converges more slowly to 0 as $\sigma\to 0$. This large deviation theoretical argument aligns with the occurrence of saddle avoidance on the extended lengthscale $\mathcal{O}(\sigma^\mu)$ (see above, \cite{maier_limiting_1997}).

Sample transition paths do not distribute proportionally to $\exp(-V_R/\sigma^2)$ across the flat quasipotential region but are skewed away from the saddle: the boundary-crossing location $\bm x_c(\varepsilon, \sigma)$ of the MTP lies at the end of $\mathcal{Z}$, far from $\bm x_M$. This is likely related to the instanton approaching the basin boundary tangentially. Yet, Kramers' scaling of the mean first-exit time $\expval{\tau_\sigma^{RL}} \asymp \exp(\sigma^{-2}\,V_R(\bm x_M))$ is empirically observed for all transition path ensembles (see Fig. S2 of the SI). We attribute this to the quasipotential barrier height at the mean crossing point $\bm x_c(\varepsilon, \sigma)$ being $V_R(\bm x_c) \approx V_R(\bm x_M) + \sigma^2$, yielding a comparable exponential scaling as long as $\sigma^2 \ll V_R(\bm x_M)$. The subexponential prefactor in Kramers' formula is approximately constant in the examined range of $\sigma$-values (not shown) and hence does not break the exponential scaling. This supports the suitability of the saddle for predicting transitions rates, even if the corresponding transition paths avoid the saddle.

\paragraph{A Finite-Noise Variational Formulation.} 
While the quasipotential can explain the occurrence of saddle avoidance, the question remains how to predict the location distribution
of transition paths for finite noise. 
If the noise is additive, there is increasing confidence \cite{thorneywork_resolution_2024, durr_onsager-machlup_1978, pinski_transition_2010, gladrow_experimental_2021} that the appropriate variational problem entails minimizing the OM action \cite{stratonovich_probability_1971, durr_onsager-machlup_1978, horsthemke_onsager-machlup_1975, pinski_transition_2010, gladrow_experimental_2021, li_gamma-limit_2021},
\begin{align} \label{eq:om}
    \tilde S_T^\sigma[\bm\varphi_t] = S_T[\bm\varphi_t] + \frac{\sigma^2}{2}\int_0^T \grad \cdot \bm b(\bm\varphi_t) \, \dd t \ ,
\end{align}
which adds a $\sigma$-dependent correction term to the FW action. This term acts as a \enquote{penalty} proportional to the divergence of $\bm b$ along the path. As we show in the SI, saddle avoidance is directly related to a positive divergence at the saddle.

\begin{figure}
    \centering
    \includegraphics[width=\columnwidth]{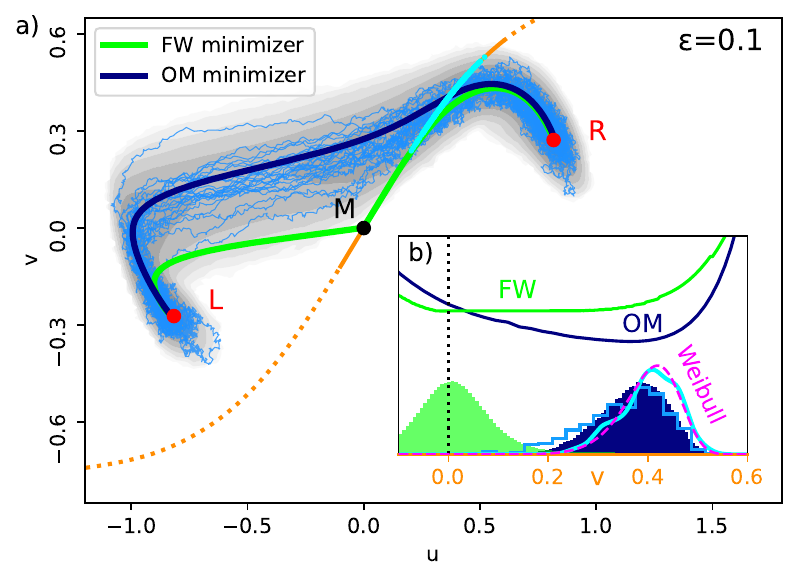}
    \caption{FHN model. a) FW (green) and OM (dark blue) minimizer compared with sample paths (light blue) for the transition $\bm\phi_t^{RL}$, showing pathspace sampling density (gray) and basin boundary (orange). b) Cross section along basin boundary (projected onto $v$). Top lines: quasipotential from minimizing FW/OM; histograms: crossing point distributions for direct sampling (light blue) vs. pathspace sampling using FW (green) and OM (dark blue) functionals; sampled first-exit point distribution (cyan) compared to rescaled Weibull distribution (magenta). }\label{fig:fw-om}
\end{figure}

Minimizing the OM action allows us to derive a candidate for the most probable transition path. Although the minimization over $T \in [0,\infty)$ is ill-defined in the OM case \cite{li_gamma-limit_2021}, one can often select a characteristic path travel time $T$; here we fix $T=\langle t_\sigma^{RL} \rangle$ to match the kinetics of the transition process. 
Numerically, we use the OM action formulation to perform \textit{pathspace sampling} (see SI), which yields a transition path density identical to that of sampling Eq. \eqref{eq:system-sde} via Monte Carlo simulation conditioned on the start and end points \cite{stuart_conditional_2004, hairer_analysis_2007}. 

We focus on the FHN model for $\varepsilon=0.1$ and $\sigma \approx 0.119$. 
The minimizer of the OM action over $\mathcal{C}_{RL}^T$ avoids the saddle and yields an appropriate approximation of the MTP (Fig. \ref{fig:fw-om}a).
Further, pathspace sampling using the OM action indeed recovers the observed transition path density (unlike using the FW action, see Fig. \ref{fig:fw-om}b), and the first-exit points are accurately described by a Weibull distribution $\rho(\sigma, \varepsilon, A)$ with $A \approx 1.23$, corroborating Ref. \cite{maier_limiting_1997}. 
The difference between the first-exit point and boundary crossing point distributions in Fig. \ref{fig:fw-om}b results from transition paths sometimes crossing the basin boundary multiple times before reaching the competing attractor.

Our results show that the OM minimizer gives better predictions of sample transition paths. This motivates constructing a finite-noise, finite-time quasipotential landscape in the spirit of Eq. \eqref{eq:fw-quasipotential} but with $T$ fixed and $\tilde S_T^\sigma$ replacing $S_T$. This quantity exhibits a minimum  on the basin boundary where sample transition paths are observed to cross (Fig. \ref{fig:fw-om}b).

\begin{figure}
    \centering
    \includegraphics[width=\columnwidth]{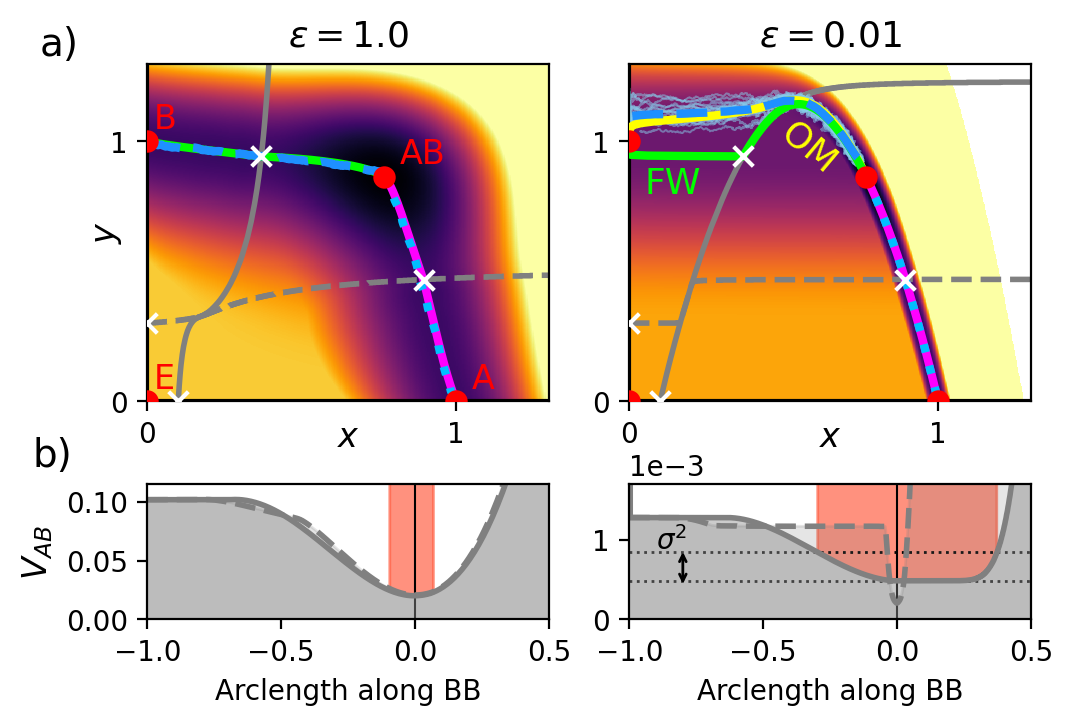}
    \caption{COMP model. a) $V_{AB}$ with respect to $\bm x_{AB}$ (brighter indicates larger $V_{AB}$) for $\varepsilon \in \{ 1, \, 0.01\}$, FW instantons to $\bm x_A$ (magenta) and $\bm x_B$ (green), and corresponding MTPs (dashed blue); OM minimizer (yellow) shown for the saddle-avoidant scenario. Basin boundaries (BB, gray lines), attractors (red points), and saddles (white crosses) are shown. b) $V_{AB}$ along the solid/dashed BB displayed in a), as a function of arclength from the respective saddle, highlighting  $\mathcal{Z}$ (red shading) for the solid BB.}
    \label{fig:comp}
\end{figure}

\paragraph{Two Competing Species.} In multiscale systems with more than two attractors, FW theory may accurately predict the path of one transition scenario while failing for another. We show this in a system of two competing species $A$ and $B$ perturbed by \textit{multiplicative} noise (COMP) \cite{bazykin_dynamics_1998}. The growth of each population is modeled with an Allee effect \cite{stephens_what_1999} according to
\begin{align}
    \bm b_{\text{comp}}(x,y) = 
    \begin{pmatrix}
    x(x-\alpha_A)(1-x) - \beta_A xy \\
    \varepsilon [ y(y-\alpha_B)(1-y) - \beta_B xy ]
    \end{pmatrix} \,.
\end{align}
Here $x$ and $y$ denote the population densities of species $A$ and $B$, respectively, which each go extinct below their critical density $\alpha_A=0.1$, $\alpha_B=0.3$. The  parameter $\varepsilon$ represents the ratio of net growth rates of the two species. The competition term is controlled by $\beta_A=0.18$ and $\beta_B=0.1$. This choice of parameters yields four stable equilibrium points: a state $\bm x_{AB}$ where both species coexist, two states $\bm x_A, \bm x_B$ where only one species survives, and a full extinction state $\bm x_E$. Additionally, four saddle points and one repeller exist in the non-negative quadrant (Fig. \ref{fig:comp}a). As with the FHN model, we investigate the dynamics of Eq. \eqref{eq:system-sde}, now for $\bm b = \bm b_\text{comp}$ and a state-dependent covariance matrix $\bm Q = \text{diag}(x, y)$ mimicking population fluctuations \cite{weissmann_simulation_2018, meng_tipping_2020}. 

Several transition scenarios are possible. We focus on the extinction of one species, realized by the two scenarios $AB \to A$ and $AB \to B$. If both species have similar net growth rates ($\varepsilon\approx 1$), both scenarios are 
characterized by a well-defined transition channel and distinct quasipotential minimum on the relevant basin boundary (Fig. \ref{fig:comp}). The MTPs track the FW instanton and cross near the corresponding saddle. Contrarily, the two scenarios differ when one species grows faster than the other (e.g. $\varepsilon=0.01$): transitions $AB \to A$ follow the instanton in a narrow quasipotential channel, whereas for $AB\to B$ the MTP detaches from the instanton and avoids the saddle in a region of flat quasipotential $V_{AB}$ induced by the fast-slow dynamics.

For the COMP system, it is possible to obtain the OM minimizer as before by performing a coordinate transform $x \to 2\sqrt{x}$, $y \to 2\sqrt{y}$ that effectively leads to an additive noise problem (see SI). Again, we find that the OM minimizer avoids the saddle and agrees with the MTP obtained from direct simulation samples (Fig. \ref{fig:comp}a).

\paragraph{Discussion.} Multiscale dynamics can cause noise-induced transitions to bypass the saddle point between competing attracting states. Since physical systems are typically non-gradient, noisy and multiscale, this phenomenon can appear in various  applications, even where transitions classify as rare events. 
Transition path ensembles may deviate strongly from the minimizer of the FW action for weak yet finite noise, while Kramers' law remains valid. These properties are possible due to a flat quasipotential along the basin boundary, which may occur due to timescale separation in the drift term for additive and multiplicative noise alike. 

Despite avoiding the instanton, sample transition paths still tend to bundle within a tube around a typical transition path, manifesting the general notion of dynamical typicality \cite{galfi_fingerprinting_2021, lucarini_typicality_2023}. For additive noise of finite amplitude, minimizing the OM action (for an appropriate path travel time $T$) yields an apt prediction of this most probable transition path and the transition path distribution, in contrast to minimizing the FW action. For multiplicative noise, computing finite-noise most probable transition paths is a topic of ongoing research \cite{thorneywork_resolution_2024, kappler_sojourn_2024, durr_onsager-machlup_1978}. Interestingly, estimating the FW minimizer with machine learing using \textit{deep gMAM} \cite{simonnet_computing_2023} yields a path that more closely resembles the OM minimizer instead of the true FW instanton \footnote{E. Simonnet, personal communication}.
Addressing the mathematical challenges of multiplicative and degenerate noise are of great interest for future work.

Ref. \cite{maier_limiting_1997} presented a path-geometric study of saddle avoidance in two dimensions based on the local stability at the saddle. Our work complements this by a) providing a global stability viewpoint based on the quasipotential, b) considering the stability of the critical manifold around the saddle, and c) clarifying the link between saddle
avoidance, positive divergence at the saddle, and multiscale dynamics. These concepts may enable anticipating saddle avoidance also in non-gradient fast-slow systems of higher dimension: if the saddle exists on an unstable branch of the critical manifold, we conjecture that repulsion away from the saddle outweighing attraction towards it will generally cause the instanton to approach the basin boundary tangentially, resulting in a quasipotential plateau along the saddle's stable set. Other physical interpretations of saddle avoidance, such as anisotropic friction \cite{berezhkovskii_rate_1989}, may be interpreted as a multiscale feature. Linking our findings to recent literature on coarse-graining \cite{Legoll2019,Hartmann2020} could provide further insights into multiscale stochastic systems.

While transitions driven by L\'evy noise are known to avoid the saddle \cite{lucarini_levy-noise_2021}, our results challenge the generic role of saddles as gateways of noise-induced transitions also under Gaussian noise. Saddle avoidance limits the classical applicability of FW theory for predicting most probable transition paths in multiscale systems, since the regime of weak but finite noise will often apply to rare events of relevance.
Especially for high-impact low-likelihood events, understanding where in state space transitions will likely occur is crucial to assess a system's resilience to random fluctuations.

\acknowledgements{The authors wish to thank P. Ditlevsen for suggesting the FHN model, M. Harsh for suggesting to transform the COMP system, P. Ashwin, M. Cameron,  K. Lux, C. Nesbitt, T. T\'el, and J. Wouters for valuable discussions, E. Simonnet for checking our results on the FHN model using deep gMAM, and the anonymous reviewers for their helpful comments. RB, RD, RR, VL and UF gratefully acknowledge funding from the European Union’s Horizon 2020 research and innovation programme under the Marie Skłodowska-Curie Grant agreement no. 956170 (CriticalEarth). VL acknowledges the support received from the EPSRC project LINK (grant no.EP/T018178/1) and from the EU Horizon 2020 project TiPES (grant no. 820970) and ClimTIP (grant no. 101137601). TG acknowledges support from EPSRC projects EP/T011866/1 and EP/V013319/1. The simulations were performed at the HPC Cluster ROSA, located at the University of Oldenburg (Germany) and funded by the DFG through its Major Research Instrumentation Programme (INST 184/225-1 FUGG) and the Ministry of Science and Culture (MWK) of the Lower Saxony State. The authors have applied a Creative Commons Attribution (CC BY) licence to any Author Accepted Manuscript version arising from this submission.}

\bibliography{saddle-avoidance}
\end{document}